\documentclass{svmult}
\usepackage{amscd,amssymb,amsmath}
\usepackage[all]{xy}

\usepackage{graphicx}
\usepackage[latin1]{inputenc}
\usepackage{longtable}
\def\C{\mathbb{C}}

\def\P{\mathbb{P}}
\def\Q{\mathbb{Q}}

\def\R{\mathbb{R}}
\def\Z{\mathbb{Z}}
\def\A{\mathbb{A}}

\def\CC{\mathcal{C}}
\def\CF{\mathcal{T}}

\def\CZ{\mathcal{Z}}
\def\Fx{F^\times}

\def\cub{\square}
\def\disjoint{\sqcup}

\def \CFT {{\widetilde{\CF}}}

\def \geq{{\,\geqslant\,}}
\def \leq{{\,\leqslant\,}}

\def \sm {{\smallskip}}

\def \om {\omega}

\def \Kriz {{K\v r\'\i\v z}}
\def \dbar{{\overline{d}}}

\newcommand{\captionfonts}{\small}
\makeatletter  % Allow the use of @ in command names
\long\def\@makecaption#1#2{%
  \vskip\abovecaptionskip
  \sbox\@tempboxa{{\captionfonts #1: #2}}%
  \ifdim \wd\@tempboxa >\hsize
    {\captionfonts #1: #2\par}
  \else
    \hbox to\hsize{\hfil\box\@tempboxa\hfil}%
  \fi
  \vskip\belowcaptionskip}
\makeatother   % Cancel the effect of \makeatletter

\newtheorem{Theorem}    {Theorem}[section]

\newtheorem{Proposition}       [Theorem]{Proposition}

\newtheorem{Lemma}      [Theorem]{Lemma}
\newtheorem{Definition} [Theorem]{Definition}

\newtheorem{Remark}[Theorem]{Remark}

\newtheorem{Example}    [Theorem]{Example}

\begin{document}
\title*{Multiple logarithms, algebraic cycles and trees}

\author{H.~Gangl\inst{1} \and A.B.~Goncharov\inst{2} \and A.~Levin\inst{3}}
\institute{MPI f\"ur Mathematik, Vivatsgasse 7, D-53111 Bonn, Germany %%\email{herbert@mpim-bonn.mpg.de} 
\and Brown University, Box 1917, Providence, RI 02912, USA 
%%\email{sasha@math.brown.edu} 
\and Institute of Oceanology, Moscow, Russia
%%\email{alevin@wave.sio.rssi.ru}
}

\maketitle

%\begin{abstract}
This is a short exposition---mostly by way of the toy models ``double logarithm'' and ``triple logarithm''---which should 
serve as an introduction to the article \cite{gangl:GGL} in which we establish a connection between multiple polylogarithms, 
rooted trees and algebraic cycles.
%\end{abstract}

% \tableofcontents
%\begin{flushleft}{\bf \small Version 1.0 (7.4.04) - DO NOT DISTRIBUTE.}
%\end{flushleft}

\section{Introduction} 

The multiple polylogarithm functions were defined in \cite{gangl:GonICM} by the power series
$$
Li_{n_1,\dots,n_m}(z_1,\dots,z_m)= \sum_{0<k_1<\cdots<k_m}
\frac{z_1^{k_1}}{k_1^{n_1}} \frac{z_2^{k_2}}{k_2^{n_2}}\dots
\frac{z_m^{k_m}}{k_{m\phantom{|}}^{n_m}}\qquad (z_i \in \C, |z_i|<1)\,.
$$
They admit an analytic continuation to a Zariski open subset of $\C^m$. 
Putting $m=1$ in this definition, we recover
 the classical polylogarithm function. 
Putting $n_1 = ... = n_m =1$, we get 
the {\it multiple logarithm} function. 

Let $x_i$ be complex numbers. Recall that an 
iterated integral is defined as %x_0\leq t_1\leq\cdots\leq t_m\leq x_{m+1}
\begin{equation}\label{int_def}
I(x_0; x_1,\dots,x_m; x_{m+1}) = \int\limits_{\Delta_\gamma} \frac{dt_1}{t_1-x_1}\wedge 
\cdots \wedge\frac{dt_m}{t_m-x_m}\,,
\end{equation}
where $\gamma$ is a path from $x_0$ to $x_{m+1}$ in $\C - \{x_1, ..., x_m\}$, and the cycle of integration  $\Delta_\gamma$ 
consists of all $m$-tuples of points $(\gamma(t_1), ..., \gamma(t_m))$ with $t_i\leq t_j$ for $i<j$.

Multiple polylogarithms can be written as iterated integrals (cf.~loc.cit.). In particular, the
 iterated integral representation of the multiple logarithm function 
 is given as 
\begin{equation} \label{ml}
Li_{1,\dots,1} (z_1,\dots,z_m)= (-1)^m I(0; x_1,\dots,x_m; 1) \,,
\end{equation}
 where we set
\begin{equation} \label{ml2}
x_1:= (z_1\cdots z_m)^{-1}, \quad x_2:= (z_2\cdots z_m)^{-1},\quad \dots\quad, x_m:= z_m^{-1}\,.
\end{equation}
Observe that in (\ref{ml2}) the parameters $x_1, \dots, x_m$ are non-zero. 
Many properties of the iterated  integrals will change if we put some of the $x_i$'s equal to zero, 
which is why the study of multiple polylogarithms cannot be directly reduced to investigating multiple logarithms only. 

{\em Notation:} We will use the notation $I_{1,\dots,1}(x_1, \dots, x_m)$ for
$I(0; x_1,\dots,x_m; 1)$ in order to emphasize that the $x_i$ are non-zero.

\medskip
In the paper \cite{gangl:BK}, Bloch and \Kriz\ defined an algebraic cycle realization 
of the classical polylogarithm function. The goal of our project \cite{gangl:GGL} was to 
develop a similar construction for multiple polylogarithms. In this paper we explain 
how to do this in the case of multiple logarithms, with special emphasis on the
toy examples of the double and triple logarithm.

The structure of the paper is as follows. Let $F$ be a field. 
In Section~2 we recall Bloch's differential graded algebra 
$\CZ^\bullet(F, \bullet)$ of cubical algebraic cycles.   
In Section~3 we define, for a set  $R$, another 
differential graded algebra $\CF_\bullet^\bullet(R)$, built from $R$-decorated rooted forests. 
A very similar differential graded algebra, for non-rooted forests, was introduced in \cite{gangl:GonArb}. 

In  Section~4 we relate these two DGA's in the case when $R = \Fx$. More precisely, we define 
a subalgebra $\widetilde \CF_\bullet^\bullet(\Fx)$ of $\CF_\bullet^\bullet(\Fx)$ 
by imposing some explicit genericity condition 
on the decoration.
Then we construct a map of graded Hopf algebras 
$$
\varphi: \widetilde \CF_\bullet^\bullet(\Fx) \longrightarrow \CZ^\bullet(F, \bullet)\,.
$$
In  Section~5 we introduce the second important ingredient of our construction:  
given a collection of elements $x_1, \dots, x_m \in \Fx$, we define an element 
$$
\tau(x_1, \ldots, x_m) \in \CF_\bullet^\bullet(\Fx)\,.
$$
Under certain explicit conditions on the $x_i$'s, it belongs to $\widetilde \CF_\bullet^\bullet(\Fx)$. 
Then the algebraic cycle $\varphi\tau(x_1, \ldots, x_m)$ corresponds to the 
multiple logarithm  (\ref{ml}). 
In Section~6, using  ideas from 
\cite{gangl:BK} concerning the Hodge realization for $\CZ^\bullet(\C,\bullet)$, 
 we show by way of example how to get the original multivalued analytic functions (\ref{ml}) %, i.e. the double and triple logarithm,
 from the constructed cycles. 

{\bf Acknowledgement}. 
This work has been done while we enjoyed the hospitality of the MPI (Bonn). We 
are grateful to the MPI for providing ideal working conditions and  
support. A.G. was supported by the NSF grant DMS-0099390. %A.L. was supported by the grant  ?????. 
A.L. was partially supported by the grant  RFFI 04-01-00642.

\section{Cubical algebraic cycles}
Let $F$ be a field. Following  \cite{gangl:Bl2}, we define the algebraic $1$-cube $\cub_F$ as a pair
$$\cub_F = \Big(\P_F^1 \setminus \{1\} \simeq \A_F^1\,, (0) - (\infty)\Big).$$
Here we consider the 
 standard coordinate $z$ on the projective line $\P_F^1$ and remove from it the point 1. Furthermore, 
 $(0) - (\infty)$ denotes 
the divisor defined by the two points $0$ and $\infty$. The algebraic $n$-cube is defined by setting 
$\cub_F^n=(\cub_F)^n$.  

Bloch defined the cycle groups
\begin{eqnarray*}
\CC^p(F,n) = \Z\big[\text{\{admissible closed irreducible subvarieties over $F$,}\\
 \text{ of codimension $p$ in $\cub_F^n$\}}\big]\,.\qquad\qquad\qquad
\end{eqnarray*}
Here a cycle is called {\bf admissible} if it intersects all the faces (of any 
codimension) of $\cub_F^n$ properly, i.e., in codimension~$p$ or not at all. 
Consider the semidirect product 
of the symmetric group $S_n$ and the group $(\Z/2\Z)^n$, acting by 
  permuting and inverting the 
coordinates in $\cub_F^n$. Let $\varepsilon_n$ be the sign representation of this group. 
The group $\CZ^p(F,n)$ is defined as the coinvariants of this group acting on  $\CC^p(F,n)\otimes \varepsilon_n$. 
Bloch showed that these groups, for a fixed ~$p$, form a complex
$$ \dots \to \CZ^p(F,n) \ { \buildrel \partial \over \to}\ \CZ^p(F,n-1) \to \dots $$
where the differential $\partial$ is given by 
$$\partial = \sum_{i=1}^n (-1)^{i-1} (\partial_0^i - \partial_\infty^i)$$
and $\partial_\varepsilon^i$ denotes 
the operator which is given by the intersection with the coordinate hyperplane 
$\{z_i=\varepsilon\}$, $\varepsilon \in \{0,\infty\}$. 

The concatenation of coordinates, followed by the corresponding 
projection under the alternation of the coordinates, provides a product on algebraic cycles,
and together with the above one gets

\begin{Proposition} (Bloch) The algebraic cycle groups $\CZ^p(F,n)$ associated to a given field $F$ provide 
a differential graded algebra $\CZ^\bullet(F,\bullet) = \sum_{p,n} \CZ^p(F,n)$.
\end{Proposition}

\noindent
\begin{Example}
For any element $a$ in $F$, one can associate such a cubical algebraic cycle corresponding to
the dilogarithm $Li_2(a)$. This cycle has been given by Totaro as 
the image of the map 
\begin{eqnarray*}
\varphi_a : \P_F^1 &\to& \big(\P_F^1\big)^3\,, \\
             t &\mapsto & (t,1-t, 1-\frac a t)\,,
\end{eqnarray*}
restricted to the algebraic cube $\cub_F^3$: we write
$$C_a := \Big[t,1-t,1-\frac a t\Big] := \varphi_a\big(\P_F^1 \big) \cap \cub_F^3\,.$$
The cycle $C_a$ belongs to the group $\CZ^2(F,3)$. One has 
$$ \partial C_a = [a,1-a]\in \cub_F^2$$
(only $\partial_0^3$ gives a non-empty contribution). 
The same computation shows that $C_a$ is in fact admissible. Observe the apparent similarity with 
the 
formula $d\, Li_2(a) = -\log(1-a)\,d\log(a)$ for the differential of  the dilogarithm. 
\end{Example}
\begin{Example}
Recall (cf. \cite{gangl:GonDL}) that the double logarithm is defined via the power series
$$Li_{1,1}(x,y) = \sum_{0<m<n} \frac{x^m}m\frac{y^n}n \qquad (|x|,|y|<1,\ x,y\in \C).$$
Its differential is computed as follows
\begin{eqnarray} \label{diffLi}
d\,Li_{1,1}(x,y) &=& \sum_{0<m<n} x^{m-1} \frac{y^n}n \, dx \ +\ \sum_{0<m<n} \frac{x^m}m y^{n-1} \, dy \nonumber \\
 &=& \sum_{n>0} \frac{1-x^{n-1}}{1-x} \frac{y^n}n \, dx \ +\ \sum_{0<m} \frac{x^m}m \, \frac{y^m dy}{1-y} \nonumber \\
 &=& Li_1(y) \frac{dx}{1-x} - Li_1(xy) \frac{dx}{x(1-x)} \ +\ Li_1(xy) \frac{dy}{1-y} \nonumber \\
 &=& Li_1(y)\, d\,Li_1(x) - Li_1(xy)\, d\,Li_1(\frac1x) + Li_1(xy)\, d\,Li_1(y)\,.
\end{eqnarray}

The cycle
\begin{equation}\label{first_ex}
C_{a,b} := \Big[1-t, 1-\frac{ab}{t},1-\frac b t\Big] \ \in \CZ^2(F,3)\,
\end{equation}
will play the role of the double logarithm $Li_{1,1}(a,b)$ among the algebraic cycles. Its boundary is readily evaluated as
\begin{equation}\label{bdry}
\partial C_{a,b} = \underbrace{\big[1-ab, 1-b\big]}_{\text{from }z_1=0} - \underbrace{\Big[1-ab, 1-\frac 1a\Big]}_{\text{from }z_2=0} + \underbrace{\big[1-b,1-a\big]}_{\text{from }z_3=0} \ \in 
\CZ^2(F,2)
\end{equation}
whose individual terms are already very reminiscent of the three terms in (\ref{diffLi}). 
\end{Example}

Observe that, setting $x_1=(ab)^{-1}$, $x_2=b^{-1}$, we can rewrite $C_{a,b}$ as
$$Z_{x_1,x_2}:=  \Big[1-\frac {1}u, 1-\frac u {x_1},1-\frac u{x_2}\Big] \,,$$
whose constants $x_1,x_2$ are chosen to match the iterated integral form $I_{1,1}(x_1,x_2)$ of $Li_{1,1}(a,b)$.
In the following, we will deal with cycles in the $Z_{\dots}$-form.  
(Note that the change of variable $t=u^{-1}$ does not change the cycle.)

Below we will explain how to generalize this definition for the case of multiple logarithms. 

\section{The differential graded algebra of $R$-deco forests}

In this paper a {\em plane tree} 
is a finite tree whose internal vertices are of valency~$\geq 3$, and where at each vertex a cyclic ordering 
of the incident edges is given.
% a cyclic ordering of the
%edges incident with a vertex.  
We assume that 
all the other vertices are of valency~1, and call them {\em external} vertices. A plane tree is {\em planted}  if it
has a distinguished external vertex of valency~1, called its {\em root}. 
A {\em forest} is a disjoint union of trees. 

\subsection{The orientation torsor}
Recall that a {\it torsor} under a group $G$ is a set on which $G$ acts freely transitively.

Let $S$ be a finite (non-empty) set. We impose on the set of orderings of $S$ an equivalence relation, given by
even permutations of the elements. The  equivalence classes form a 2-element set ${\rm Or}_S$. It 
has an obvious $\Z/2\Z$-torsor structure and is called the {\it orientation
torsor of $S$}. 

\begin{Definition} 
Let $F$ be a plane forest. The {\bf orientation torsor of} $F$ is the orientation
torsor of the set of its edges. 
\end{Definition}

Observe that once we have chosen an edge ordering in a plane tree $T$, 
e.g., by fixing a root, there is a canonical orientation on $T$.

\subsection{The algebra of $R$-deco forests}
\begin{Definition}\label{Rdeco} Let $R$ be a set. An {\bf $R$-deco tree} is a planted plane tree with a map, 
called {\bf $R$-decoration}, from its external vertices to $R$. An {\bf $R$-deco forest} is a disjoint 
union of $R$-deco trees.  
\end{Definition}
\begin{Remark}\label{rem1}
\begin{enumerate}
\item There is an obvious induced direction for each edge in an $R$-deco tree, away from the root.
\item There is an ordering of the edges, starting from the root edge, which is 
induced by the cyclic ordering of edges at internal vertices.
\end{enumerate}
\end{Remark}
Our convention for drawing the trees is that the cyclic ordering of edges around
internal vertices is displayed in counterclockwise direction.

\begin{Example} We draw an $R$-deco tree $T$ with root vertex decorated by $x_4\in R$; its other external vertices are decorated
by $x_1$, $x_2$, $x_3\in R$. The above-mentioned ordering of the edges $e_i$ coincides
with the natural ordering of their indices, while the direction
of the edges (away from the root) is indicated by small arrows along the edges.
\end{Example}
$$
\includegraphics{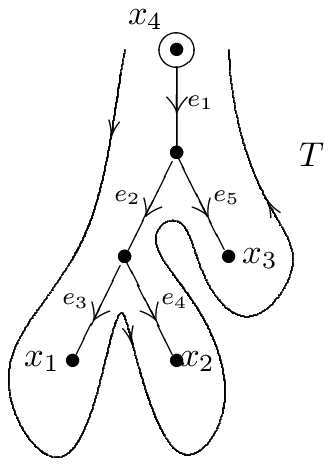}
$$
%%%\begin{figure} \label{treeart} \hskip 0pt
%\vskip 15pt
%\hskip 100pt 
%\xy  0;<-3pt,0pt>:
%(5,0)*+{};(-5,-1)*+{}
%%% **\crv{(10,20)&(15,30)&(10,-10)&(-10,-20)}
 %**\crv{(6,5)&(7,20)&(20,30)&(10,45)&(6,25)&(5,25)&(4,30)&(0,42)&(-6,35)&(-3,25)&(4,18)&(-1,15)&(-2,22)&(-7,28)&(-13,20)&(-8,15)&(-5,5)}
%?(.03)*\dir{>} ?(.35)*\dir{>} ?(.93)*\dir{>}
 %\POS(0,0) *+{\bullet} *\cir{} 
 %\ar @{-} +(0,10)*{\bullet}^{e_1} \POS(0,6)*\dir3{>}
%%% ?(.5)*\dir{>}
  %\POS(0,10) *{\bullet}
 %\ar @{-} +(5,10)*{\bullet}_{e_2} \POS(3,16)*\dir3{>}\POS(0,10)
 %\ar @{-} +(-5,10)*{\bullet}^{e_5} \POS(-3,16)*\dir3{>}
 %\POS(5,20) *{\bullet}
 %\ar @{-} +(5,10)*{\bullet}_{e_3} \POS(8,26)*\dir3{>}\POS(5,20)
 %\ar @{-} +(-5,10)*{\bullet}^{e_4} \POS(2,26)*\dir3{>}
%\POS(3,-3) *{x_4}
%\POS(13,30) *{x_1}
%\POS(-2,30) *{x_2}
%\POS(-8,20) *{x_3}
%%%\POS(12,15) *{\gamma}
%\POS(-13,10) *{T}
%%% **\crv{(10,50)&(30,20)&(50,-20)&(60,-10)}
%\endxy
%\caption{

\medskip\noindent
%%%\caption
\centerline{\small Figure 1: An $R$-deco tree $T$ with root vertex decorated by $x_4$.}
%% ; its other external vertices are decorated by $x_1$, $x_2$, $x_3$.
%%% \end{figure}
\vskip 15pt

\begin{Definition} 
For a set $R$, the $\Q$-vector space ${\CF}_{\bullet}(R)$ is generated by the elements $(F, \omega)$, where $F$ is an $R$-deco forest 
and $\omega$ an orientation on it, subject to the relation $(F, -\omega) = -(F, \omega)$. 
\end{Definition}

We define the {\it grading} of an $R$-deco tree $T$  by
$$e(T) = \#\{\text{edges of }T\}\,$$ 
and extend it to forests by linearity: $e(F_1\disjoint F_2): = e(F_1) + e(F_2)$. 
It provides the vector space ${\CF}_{\bullet}(R)$ with a natural grading. 

We define the algebra structure $\star$ on ${\CF}_{\bullet}(R)$ by setting
$$
(F_1, \omega_1)\star (F_2, \omega_2) := (F_1\disjoint F_2, \omega_1 \otimes \omega_2)
$$
It makes ${\CF}_{\bullet}(R)$ into a graded commutative algebra,  
called {\it the $R$-deco forest algebra $\CF_{\bullet}(R)$}.

Let $V_{\bullet}({R})$ be the graded $\Q$-vector space with  basis given by $R$-deco trees, with  
the above grading. 

\begin{Lemma} \label{fga}
The  algebra $\CF_{\bullet}(R)$ is the free graded commutative algebra 
generated by the graded vector space $V_{\bullet}({R})$.
\end{Lemma}

So the basis elements of $V_{\bullet}({R})$ commute in the $R$-deco forest algebra via the rule
$$ (T_1,\omega_{1})\star (T_2,\omega_{2}) = (-1)^{e(T_1)e(T_2)} (T_2,\omega_{2})\star (T_1,\omega_{1})\,.$$

\subsection{The differential}
A differential on ${\CF}_{\bullet}(R)$ is a map
$$
d: {\CF}_{\bullet}(R) \longrightarrow {\CF}_{\bullet-1}(R) 
$$
satisfying $d^2=0$ and the Leibniz rule. Since, by Lemma \ref{fga}, ${\CF}_{\bullet}(R)$ is a free graded commutative algebra, 
it is sufficient to define it on the algebra generators, that is on the elements $(T, \omega)$, where $T$ is an $R$-deco tree
and $\omega$ is an orientation of $T$.

The terms in the differential of a tree $T$ arise by contracting an edge of~$T$---they fall into two types,
according to whether the edge is internal or external. We will need the notion of a splitting.

\begin{Definition} A {\bf splitting} of a tree $T$ at an internal vertex $v$ is the disjoint union of 
the trees which arise as $T_i\cup v$ where the $T_i$ are the connected components of $T\setminus v$.

Further structures on $T$, e.g.~a decoration at $v$, planarity of $T$ or an ordering of its edges, are inherited for each $T_i\cup v$.
Also, if $T$ has a root $r$, $v$ plays the role of the root for all $T_i\cup v$ which do not contain $r$.
\end{Definition} 

\begin{Definition} Let $e$ be an edge of a tree $T$. The {\bf contraction} of $T$
along $e$, denoted $T/e$, is given as follows:
\begin{enumerate}
\item If $e$ is an internal edge, then $T/e$ is again a tree: it is the same tree as $T$ 
except that $e$ is contracted and the incident vertices $v$ and $v'$ of $e$ are identified 
to a single vertex.
\item If $e$ is an external edge, then $T/e$ is obtained as follows: first we contract
the edge $e$ to a vertex $w$ and then we perform a splitting at $w$.
\end{enumerate}
\end{Definition} 

Two typical examples are given below: in %% the pictures below,
Figure~2 we contract a {\em leaf}, i.e.~an external vertex which is not the root vertex,
and in Figure~3 the root vertex is contracted.
%%\begin{figure}

$$
\includegraphics{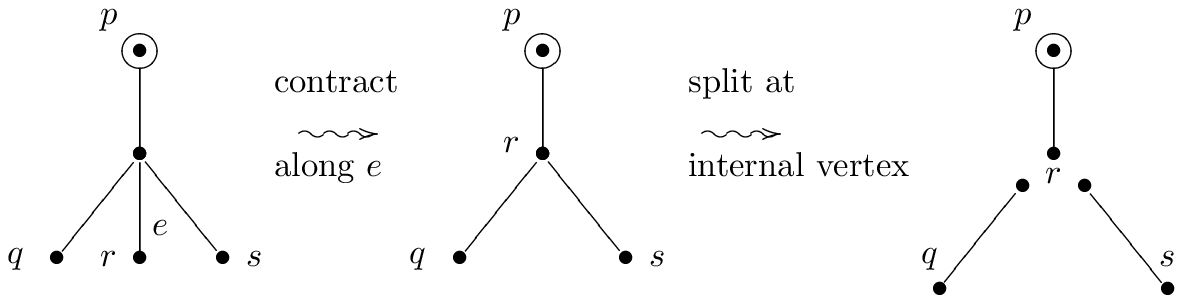}
$$

%%%\caption
\centerline{\small Figure 2: Contracting a leaf.}
%%\end{figure}

%%%\begin{figure}
$$
\includegraphics{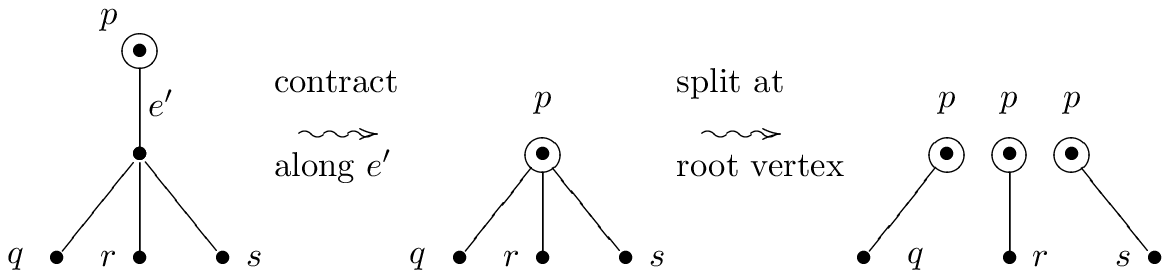}
$$
%%%\caption
\centerline{\small Figure 3: Contracting the root edge.}
%%%\end{figure}
\vskip 15pt

Let $S$ be a finite set, and ${\rm Or}_S$ the orientation torsor of $S$. 
We present its elements as $s_1 \wedge ... \wedge s_n$, where $n = |S|$. Now given an element $s \in S$ and 
$\omega\in {\rm Or}_S$, we define 
an element $i_s\omega \in {\rm Or}_{S-s}$ as follows:
$$
i_s\omega:= s_2 \wedge ... \wedge s_n \quad \mbox{if $\omega =  s \wedge s_2 \wedge ... \wedge s_n$}\,.
$$

\begin{Definition} Let $T$ be a finite tree with set of edges $E$, and let $\om$ be an orientation of $T$. 
The {\bf differential} on $(T,\om)$ is defined as
$$d: (T,\om) \mapsto \sum_{e\in E}(T/e, i_e\om)\,.$$
%%If, furthermore, $T$ is an $R$-deco tree, then the 

\end{Definition}

\begin{Example}
The simplest non-trivial example for the differential of an $\Fx$-deco tree, where $\Fx$ is the multiplicative group of a field $F$, can be seen on a tree with one internal vertex, 
as given in Figure~4. Here we choose the $\Fx$-decoration $(x_1,x_2)$ with $x_1,x_2\in \Fx$ for the leaves and the
decoration $1$ for the root:
\vskip 3pt
%%%\begin{figure}
$$
\includegraphics{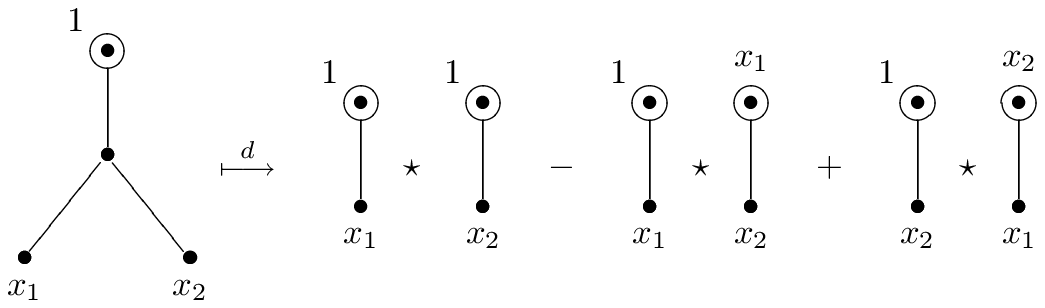}
$$
%%%\caption
\centerline{\small Figure 4: The differential on a tree with one internal vertex.}
%%%\end{figure}

%A comparison with (\ref{bdry}) gives us a first idea how to relate the trees to algebraic cycles, and this is explained
%in the next section.
\vskip 7pt
\noindent For the drawings, we use the canonical ordering of edges for $\Fx$-deco 
trees and the induced
ordering for forests which arise from a splitting. The ordering of the forest 
encodes its orientation, i.e. the choice of an element in the orientation torsor. 
\end{Example}

\begin{Remark}
There is a $\Z$-bigrading on $\CF_{\bullet}(R)$ given by the groups
$$
\CF_n^p(R) = \Z[\{\text{$R$-deco forests with $n$ edges and $p$ leaves}\}]\,.
$$
It will correspond below to the bigrading on the cycle groups $\CZ^p(F,n)$. 
\end{Remark}

For a set $R$, put
$$
\CF_{\bullet}^\bullet(R):= \bigoplus_{p\geq 0}\ \bigoplus_{0\leq n\leq p}\CF_n^p(R)\,.
$$
With the above definitions, we have:
\begin{Proposition} 
$(\CF_{\bullet}^\bullet(R), d)$ is a bigraded, differential graded algebra.
\end{Proposition}

%$[$Drop orientation $\omega$ from the pair $(F,\omega)$, once the forests have a canonical orientation,
%e.g. by defining an order on the set of its trees.$]$ 

\section{Mapping forests to algebraic cycles}

In the special case where $R=\Fx$, the multiplicative group of a field $F$, we can 
establish the connection between the two differential graded algebras
above, and Theorem \ref{dga_morphism} below gives the main result of this paper.

It turns out that the admissibility condition on algebraic cycles mentioned above forces us to restrict
to a subalgebra of $\CF^\bullet$, which we now describe.

\begin{Definition}
We call an $R$-deco tree  {\bf generic} if all the individual 
decorations of external vertices are different.

We denote the subalgebra of $\CF^\bullet_\bullet(R)$ generated by generic $R$-deco trees by $\CFT^\bullet_\bullet(R)$.
\end{Definition}

One of our key results is the following statement:
\begin{Theorem}\label{dga_morphism} For a field $F$, there is a natural map of differential graded algebras
$$\CFT^\bullet_\bullet(\Fx) \to \CZ^\bullet(F,\bullet)\,.$$
It is given by the map in Definition \ref{fc} below.
\end{Theorem}

\begin{Definition} \label{fc}
The {\bf forest cycling map} for a field $F$ is the map $\varphi$ from $\CF^\bullet_\bullet(\Fx)$ to (not necessarily admissible) 
cubical algebraic cycles over $F$ given on generators, i.e. $\Fx$-deco trees, as follows:
\begin{enumerate}
\item to each internal vertex $v$ of $T$ we associate a decoration consisting of an independent (``parametrizing'') variable;
\item to each edge with (internal or external) vertices $v$ and $w$ equipped with respective decorations $y_v$ and $y_w$
(variables or constants) we associate the expression $[1-y_w/y_v]$ as a parametrized coordinate in $\P_F^1$;
\item choosing an ordering of edges of $T$ corresponding to $\om$, we concatenate all the respective coordinates 
produced in the previous step.
%%$[$TO BE CHANGED$]$ in the canonical edge ordering of $T$, we  concatenate all the coordinates produced in the previous step.
\end{enumerate}
\end{Definition}

This somewhat lengthy description is easily understood by looking at an example. 
We will denote the concatenation product for algebraic cycles by~$*$, and
we encode the expression $1-\frac x y$, for $x,y$ in a field, by the following picture
\vskip 5pt 
$$
\includegraphics{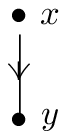}
$$

\begin{Example} Let us consider the forest cycling map $\varphi$ for the following $R$-deco tree $(T,\om)$, where the 
orientation $\om$ is given by $e_1\wedge e_2\wedge e_3$ (we leave out the arrows 
since the edges are understood to be directed away from the root):
\vskip 10pt
%%%\begin{figure}
$$
\includegraphics{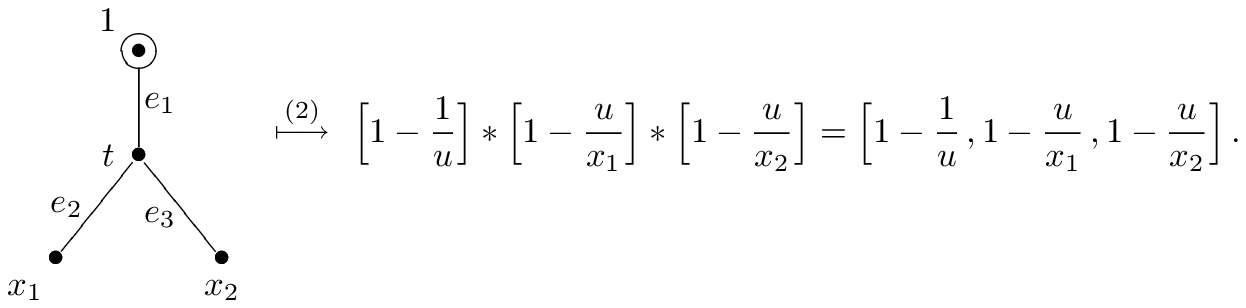}
$$
\centerline{\small Figure 5: The forest cycling map on a tree with one internal vertex.}
%%%\end{figure}
 \vskip 10pt\noindent
\end{Example}
\noindent
This cycle, as already mentioned, corresponds to the double logarithm $I_{1,1}(x_1,x_2)$, as we will see in Section~6.

\begin{Lemma} Each generic $\Fx$-deco tree maps to an admissible cycle in $\CZ^\bullet(F,\bullet)$.
\end{Lemma}
For a proof, we refer to \cite{gangl:GGL}; the main idea is that at each internal vertex we have at least one incoming
and one outgoing edge, which implies that their respective coordinates in the associated cycle cover up for each other.

\section{The algebraic cycle corresponding to the multiple logarithm}

\begin{Definition} \label{ml11}
Let $\{x_1, \dots, x_m\}$ be a collection of distinct elements of $\Fx\setminus\{1\}$. 
Then $\tau(x_1, \dots, x_m)$ is the sum of all trivalent $\Fx$-deco trees with $m$ leaves whose $\Fx$-decoration is
given by $(x_1, x_2,\dots, x_{m-1}, x_m)$, while the root is decorated by~$1$.  
\end{Definition} 
%$Li_{1,\dots,1}(z_1,\dots,z_n)$ 
Recall that the number of such trees is given by the Catalan number $\frac 1{m}{2(m-1) \choose m-1}$.

\vskip 8pt\noindent
Combining Definition \ref{ml11} with the forest cycling map $\varphi$, we get the cycles 
corresponding to the multiple logarithms. For $m=2$, the tree $\tau(x_1,x_2)$ is given by the tree
on the left of Figure~5, where $x_3$ is put equal to 1.

%\begin{Example} Let $m=2$. The cycle $C_{a,b}\in \CZ^2(F,3)$ for the
%double logarithm $Li_{1,1}(a,b)=I\big(0;(ab)^{-1},b^{-1};1\big)$ was given in (\ref{first_ex}). 
%The corresponding tree $\tau\big((ab)^{-1},b^{-1}\big)$ is given by the 
%tree on the left of Figure~4.  
%\end{Example}

\begin{Example}
The simplest example where the sum of trees consists of more than a single term 
appears when $m=3$:

\vskip 10pt
%%%\begin{figure}

$$
\includegraphics{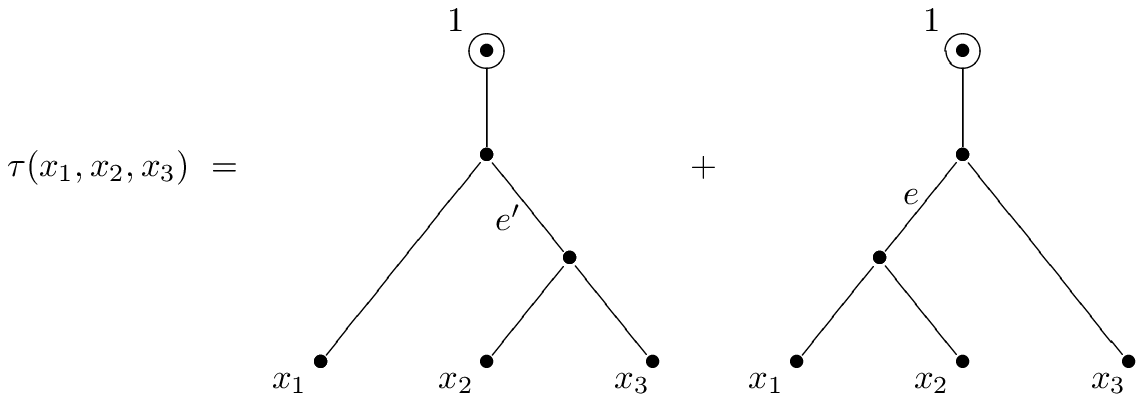}
$$
\centerline{\small Figure 6: The sum of trees corresponding to the triple logarithm $I_{1,1,1}(x_1, x_2, x_3)$.}
%%%\end{figure}
\vskip 8pt
\end{Example}

Applying $\varphi$, we get the following
cycle corresponding to the triple logarithm $I_{1,1,1}(x_1,x_2,x_3)$: 
$$  %%\begin{equation*}
Z_{x_1,x_2,x_3} = \Big[1-\frac 1 t,1-\frac t{x_1},1-\frac t u,1-\frac u{x_2} ,1-\frac u {x_3}\Big]
\ + \ \Big[1-\frac 1 t, 1-\frac t u, 1-\frac u{x_1}, 1-\frac u{x_2} ,1-\frac t {x_3}\Big]. \,
$$  %%\end{equation*}
%%%CHANGE TO Z
Here we have two parametrizing variables, $t$ and $u$, and $Z_{x_1,x_2,x_3}\in \CZ^3(F,5)$.

\begin{Example}
In Figure~7, we give the sum of $\Fx$-deco trees associated to the weight~4 multiple logarithm 
$I_{1,1,1,1}(x_1,x_2,x_3,x_4)$. %It is a formal sum of 5 trees.
%% \vskip 50pt
\vskip 15pt
%%\hskip 100pt 
%%%\begin{figure}
$$
\includegraphics{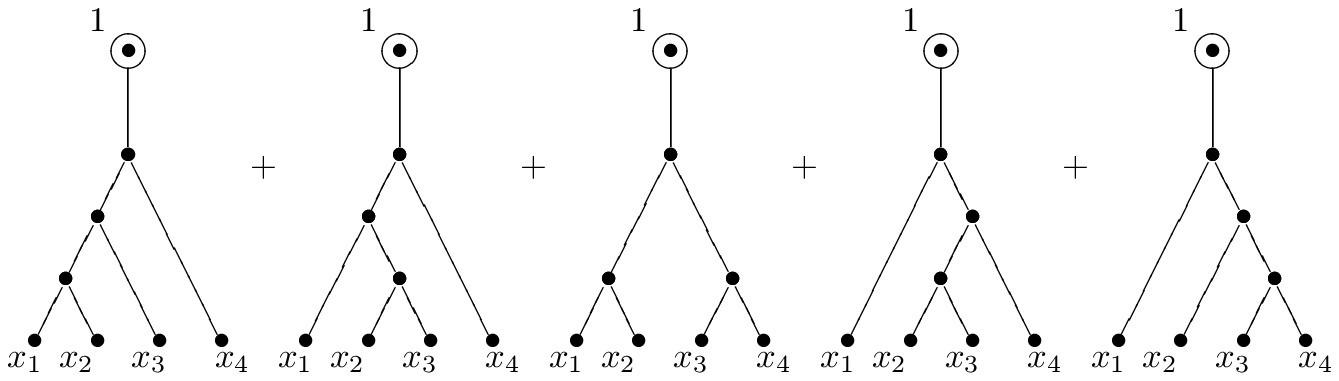}
$$
\centerline{\small Figure 7: The sum of $5$ trees corresponding to the multiple 4-logarithm.}
%%% $Li_{1,1,1,1}(z_1,z_2,z_3,z_4)$.}
%%%\end{figure}
\end{Example}
\vskip 8pt

One of the crucial properties of the elements 
 $\tau(x_1, ..., x_m)$ is that the contributions of internal edges to the differential of 
$\tau(x_1, ..., x_m)$ cancel pairwise.  This  property 
ensures, cf.~\cite{gangl:GGL}, that one can build from 
$\tau(x_1, ..., x_m)$ an element in the
Hopf algebra $\chi_{\rm mot}$ of~\cite{gangl:BK}.
Here is how the cancellations take place for $m=3$. The  
terms in $d\tau({x_1,x_2,x_3})$ coming from contracting
the {\em internal} edges $e$ and $e'$ cancel each other: 
$e$ is the third one in the (canonical counterclockwise) ordering of edges
in the first tree, while $e'$ is the second edge in the second tree, 
so the corresponding signs in the differential are
opposite. 
This means that each term of 
 $d\tau({x_1,x_2,x_3})$ is decomposable (in fact, since 
$\tau({x_1,x_2,x_3})$ consists of trivalent trees, each such term is a 
product of precisely two trees).

\section{Associating integrals to the multiple logarithm cycles}
So far we have not given a reason why we can consider the cycles associated to certain trees as ``avatars'' of 
multiple logarithms. 
In this section, we indicate how we can associate an integral to the cycles $Z_{x_1,x_2}$ and  $Z_{x_1,x_2,x_3}$ 
which is nothing else than the integral
presentation for the double and triple logarithm, respectively. 

Following Bloch and \Kriz, we embed the algebraic cycles into a larger 
set-up of ``hybrid'' cycles which have both algebraic
and topological coordinates as well as both types of differentials, and then apply the bar construction. 
We only consider ``topological'' variables $s_i\in [0,1]\subset \R$ 
subject to the condition $s_i\leq s_j$ if $i<j$,
%%0=s_0\leq s_1\leq\dots\leq s_r\leq s_{r+1}=1$, 
and taking the topological boundary~$\delta$ for a cycle with topological dimension $r$ 
amounts to taking the formal alternating sum over the 
subvarieties where either $s_k=s_{k+1}$ for some $k=1,\dots,r-1$ or $s_1=0$ or $s_r=1$.

\begin{Example} \begin{enumerate}
\item 
In order to bound $Z_{x_1,x_2}$, consider the algebraic-topological cycle parametrized by 
$t\in \P_F^1$ and $s_1\in \R$,  $0\leq s_1\leq 1$, as
$$\Big[1-\frac {s_1}t,1-\frac t{x_1},1-\frac t{x_2}\Big]\,,$$
whose topological boundary terms are obtained by putting $s_1=0$ (which produces the empty cycle) or by $s_1=1$
which yields $Z_{x_1,x_2}$. Its algebraic boundary is given by 
$$[1- \frac{s_1}{x_1},1-\frac {s_1}{x_2}] - [1-\frac{s_1}{x_1},1-\frac {x_1}{x_2}] + [1-\frac{ s_1}{x_2},1-\frac{x_2}{x_1}]\,,\quad 0\leq s_1\leq 1\,,$$
where the last two terms are ``negligible'' for the following.
\item Consider the topological cycle parametrized by $0\leq s_1\leq s_2\leq 1$, $s_i\in \R$, as
$$[1-\frac{s_1}{x_1}, 1-\frac{s_2}{x_2}]\,,$$
whose boundary terms arise from setting $s_1=0$, $s_1=s_2$ or $s_2=1$, giving the empty cycle, %% $\emptyset$, 
$[1-\frac{s_1}{x_1},1-\frac{s_1}{x_2}]$ or $[1-\frac{s_1}{x_1}, 1-\frac 1{x_2}]$, respectively.
\end{enumerate}
\end{Example}
What we are after is a cycle $\eta$ in this larger (algebraic-topological) cycle complex which bounds the cycle $Z_{x_1,x_2}$
i.e.,~such that $Z_{x_1,x_2}=(\partial+\delta)\eta$, since the ``bounding process'' will give rise to a purely topological
cycle against which we can then integrate the standard volume form $(2\pi i)^{-r}\displaystyle \frac {dz_1}{z_1}\wedge\dots\wedge \frac{dz_r}{z_r}$
(here $r=2$), $z_i$ being a coordinate on the $i^{\rm th}$ cube $\cub_\R$.

In fact, working modulo the ``negligible'' term $[1-\frac{ s_1}{x_1},1-\frac{1}{x_2}]$ above we get that the full differential $D=\partial+\delta$
on the algebraic-topological cycle groups gives 
\begin{equation}\label{boundCab}
Z_{x_1,x_2} = D\Big( \big[1- \frac{s_1}t,1-\frac t{x_1},1-\frac t{x_2}\big] + [1-\frac{s_1}{x_1},1-\frac{s_2}{x_2}]\Big)
\end{equation}
(two of the boundary terms cancel) and we associate to $Z_{x_1,x_2}$ the integral
\begin{multline}
\frac 1{(2\pi i)^2}\int\limits_{[1-\frac{s_1}{x_1},1-\frac{ s_2}{x_2}]\atop 0\leq s_1\leq s_2\leq 1}{} \frac {dz_1}{z_1}\wedge \frac{dz_2}{z_2} = 
\frac 1{(2\pi i)^2}\int\limits_{0\leq s_1\leq s_2\leq 1} \frac {d(1-\frac{ s_1}{x_1})}{1-\frac{ s_1}{x_1}}
\wedge \frac{d(1-\frac{s_2}{x_2})}{1-\frac{ s_2}{x_2}} \\ = \frac 1{(2\pi i)^2}\int\limits_{0\leq s_1\leq s_2\leq 1} \frac{ds_1}{s_1-x_1} \wedge \frac{ds_2}{s_2-x_2}\,.
\end{multline}
Therefore we see that the cycle $Z_{x_1,x_2}$ corresponds---in a rather precise way---to the iterated
integral $I_{1,1}(x_1,x_2)$.
\sm

Note that the three ``negligible'' terms encountered above can also be covered as part of a boundary if we introduce,
following \cite{gangl:BK}, yet another differential $\dbar$ (coming from the well-known bar construction), and in the ensuing
tricomplex all the terms above are taken care of. With the usual bar notation $\big|$ for a certain tensor product, 
the correct cycle is given by
\begin{multline*}
[1-\frac{s_1}{t},1-\frac t{x_1},1-\frac t{x_2}] + [1-\frac {s_1}{x_1},1-\frac{s_2}{x_2}] \\
+\ \Big(\,[1-\frac {s_1}{x_1}\big| 1-\frac{1}{x_2}] - [1-\frac {s_1}{x_1}\big| 1-\frac{x_2}{x_1}] + [1-\frac {s_1}{x_2}\big| 1-\frac{x_1}{x_2}]\,\Big)\,,
\end{multline*}
and its image under the boundary $\partial + \delta+\dbar$ is precisely the ``bar version'' of $-Z_{x_1,x_2}$, given by 
$-Z_{x_1,x_2} + \big(\,[1-\frac{1}{x_1}\mid 1-\frac{1}{x_2}] - [1-\frac{1}{x_1}\mid 1-\frac{x_1}{x_2}] + [1-\frac{1}{x_2}\mid 1-\frac{x_2}{x_1}]\,\big)$. For details,
we refer to \cite{gangl:GGL}.

\subsection{The integral associated to the triple logarithm}
In a similar fashion as for $Z_{x_1,x_2}$, we can find a bounding cycle $\gamma$ in the larger
algebraic-topological groups. We picture the see-saw-like process with the main terms (all terms which do not appear
in the diagram are ``topologically decomposable'' and therefore negligible for the final integral).
In the cycles below, we have $t,u\in\P_F^1,$  $(t,u)\not=(0,0)$, and the range of the parameters $0\leq s_i\leq 1$ is given
via the inequalities $s_i\leq s_j$ for $i<j$.

%%\begin{figure}
$$\begin{array}{rll}
%%\\
\phantom{\Bigg|\Bigg|}{\displaystyle \Big[1-\frac 1 t,1-\frac t{x_1},1-\frac t u,1-\frac u{x_2} ,1-\frac u {x_3}\Big]} & &\\
{\displaystyle + \Big[1-\frac 1t, 1-\frac t u, 1-\frac u{x_1}, 1-\frac  u{x_2} ,1-\frac t{x_3}\Big]} & &\\
 & \hskip-5pt\nwarrow^{\hskip-5pt-\delta}& \\
&{\displaystyle \phantom{\Bigg|\Bigg|}\Big[1-\frac{s_1}t,1-\frac t{x_1},1-\frac t u,1-\frac  u{x_2} ,1-\frac u{x_3}\Big]}\\
&{\displaystyle  + \Big[1-\frac{s_1}{t}, 1-\frac t u, 1-\frac u{x_1}, 1-\frac  u{x_2} ,1-\frac t{x_3}\Big]}\\
 & \hskip-5pt\swarrow_\partial& \\
{\displaystyle \phantom{\Bigg|\Bigg|}\Big[1-\frac{s_1}{x_1},1- \frac{s_1}{u},1-\frac  u{x_2} ,1-\frac u{x_3}\Big]}&\\
{\displaystyle  + \Big[1-\frac {s_1}u, 1-\frac u{x_1}, 1-\frac  u{x_2} ,1- \frac{ s_1}{x_3}\Big]}&\\
 & \hskip-5pt\nwarrow^{\hskip-5pt-\delta}& \\
 &{\displaystyle \phantom{\Bigg|\Bigg|}\Big[1-\frac{s_1}{x_1},1- \frac{s_3}{u},1-\frac  u{x_2} ,1-\frac u{x_3}\Big]}\\
 &{\displaystyle  + \Big[1-\frac {s_1}u, 1-\frac u{x_1}, 1-\frac  u{x_2} ,1- \frac{ s_3}{x_3}\Big]}\\
 & \hskip-5pt\swarrow_\partial& \\
{\displaystyle \phantom{\bigg|}-\Big[1-\frac{s_1}{x_1},1- \frac{s_3}{x_2} ,1-\frac{s_3}{x_3}\Big]}&\\
{\displaystyle \phantom{\bigg|} + \Big[1-\frac{s_1}{x_1},1- \frac{s_1}{x_2} ,1-\frac{s_3}{x_3}\Big]}&\\
 & \hskip-5pt\nwarrow^{\hskip-5pt-\delta}& \\
&{\displaystyle \phantom{\bigg|}\Big[1-\frac{s_1}{x_1}, 1-\frac{s_2}{x_2}  ,1-\frac{s_3}{x_3}\Big]}\\
\end{array}
$$
%%\vskip 5pt
\centerline{\small Figure 8: A see-saw procedure to find the integral for the triple logarithm cycle.}
%%\centerline{\small  The real variables $s_i\in [0,1]$ satisfy $s_i\leq s_j$ for $i<j$.}
\vskip 10pt
%Whenever the real variables $s_i$ appear in a cycle,
%they are understood to lie between~0 and~1, and to further satisfy $s_i\leq s_j$ for $i<j$.}
%%\end{figure}

The explanation of this diagram is as follows: the first two lines give the algebraic cycles for $Z(x_1,x_2,x_3)$; they are the images
(modulo decomposable cycles) under $\delta$ of the algebraic-topological cycles shown in the next two lines (a real parameter 
$s_1$ enters). If we apply the full differential $D=\partial+\delta$ to the latter two cycles we obtain two new irreducible cycles, 
listed in the following two lines. Again we can bound the latter (under $\delta$) by two algebraic-topological cycles 
(a new parameter $s_3$ enters, and we have $0\leq s_1\leq s_3\leq 1$; note that only the $\delta$-boundaries from putting $s_1=s_3$ are
indecomposable). In the same fashion we consider the $\partial$-boundaries of the latter, winding up with 
$\gamma_{1,3}=-[1-\frac{s_1}{x_1},1-\frac{s_3}{x_2} ,1-\frac{s_3}{x_3}]+ [1-\frac{s_1}{x_1}, 1-\frac{s_1}{x_2} ,1- \frac{s_3}{x_3}]$, $0\leq s_1\leq s_3\leq 1$ (note that the
signs are opposite). In the final step, we see that the indecomposable part of the $\delta$-boundary of the purely
topological cycle
$[1-\frac{s_1}{x_1}, 1- \frac{s_2}{x_2} ,1- \frac{s_3}{x_3}]$ $(0\leq s_1\leq s_2\leq s_3\leq 1)$ is precisely the above cycle $\gamma_{1,3}$.
Therefore the sum over the five cycles on the right hand side of the picture bounds $-Z_{x_1,x_2,x_3}$, and the only 
cycle which gives a non-trivial integral against the standard volume form $(2\pi i)^{-3}\frac{dz_1}{z_1}\wedge 
\frac{dz_2}{z_2}\wedge\frac{dz_3}{z_3}$ is the purely topological one, providing the integral
$$-\frac 1{(2\pi i)^3} \int\limits_{0\leq s_1\leq s_2\leq s_3\leq 1} \frac{ds_1}{s_1-x_1} \wedge \frac{ds_2}{s_2-x_2} \wedge \frac{ds_3}{s_3-x_3}\,,$$
which is nothing else---up to the normalizing factor $(2\pi i)^{-3}$---than $I_{1,1,1}(x_1,x_2,x_3)$.

\section{Outlook}

We can also describe multiple {\em poly}logarithms $Li_{k_1,\dots,k_n}(z_1,\dots,z_n)$ in a similar fashion, but we need to 
introduce trees with two different
kinds of external edges, and an accordingly modified forest cycling map gives us the associated algebraic cycles (cf.~\cite{gangl:GGL}).

A general construction in \cite{gangl:BK} based on the well-known bar construction applies to both DGA's, providing in particular a Hopf algebra
structure on the $R$-deco forests, and $\varphi$ furthermore induces a morphism of Hopf algebras. Details for this, as well
as the connection to the ``motivic world'', are given in \cite{gangl:GGL}.


\begin{thebibliography}{99}

%%\bibitem{gangl:BlochAlgCyc}  {\bf  Bloch, S.} {\it Algebraic cycles and higher $K$-theory.}  Adv. in Math.  61  (1986),  no. 3, 267--304. 

\bibitem{gangl:BK} {\bf  Bloch, S., \Kriz, I.} {\it Mixed Tate motives.} Ann. of Math. (2)  140  (1994),  no. 3, 557--605. 

%%\bibitem{gangl:Bl1} {\bf Bloch, S.} {\it Algebraic cycles and the Lie algebra of mixed Tate motives.}  J. Amer. Math. Soc.  4  (1991),  no. 4, 771--791. 

\bibitem{gangl:Bl2} {\bf Bloch, S.} {\it  Lectures on mixed motives.}  Algebraic geometry---Santa Cruz 1995,  329--359, Proc. Sympos. Pure Math., 62, Part 1, Amer. Math. Soc., Providence, RI, 1997.

%%\bibitem{gangl:DG} {\bf Deligne, P., Goncharov, A.B.} {\it Groupes fondamentaux motiviques de Tate mixte.} arXiv: math.NT/0302267.

\bibitem{gangl:GGL} {\bf Gangl, H.; Goncharov, A.B.; Levin, A.} {\it Multiple polylogarithms, polygons, trees and algebraic cycles.} In preparation.

\bibitem{gangl:GMS} {\bf Gangl, H.; M\"uller-Stach, S.} {\it Polylogarithmic identities in cubical higher Chow groups.}  Algebraic $K$-theory (Seattle, WA, 1997),  25--40, Proc. Sympos. Pure Math., 67, Amer. Math. Soc., Providence, RI, 1999. 

\bibitem{gangl:GonICM} {\bf Goncharov, A.B.} {\it Polylogarithms in arithmetic and geometry}. Proceedings of the International
  Congress of Mathematicians, Vol. 1, 2 (Zürich, 1994), 374--387, Birkhäuser, Basel, 1995. 

\bibitem{gangl:GonDL} {\bf Goncharov, A.B.}  {\it The double logarithm and Manin's complex for modular curves.}  Math. Res. Lett.  4  (1997),  no. 5, 617--636. 

\bibitem{gangl:GonArb} {\bf Goncharov, A.B.}  {\it Galois groups, geometry of modular varieties and graphs.} Talk at the 
Arbeitstagung Bonn, 1999. Available under  http://www.mpim-bonn.mpg.de/html/preprints/preprints.html

\bibitem{gangl:GonBar} {\bf Goncharov, A.B.}  {\it Multiple $\zeta$-values, Galois groups, and geometry of modular varieties.}  European Congress of Mathematics, Vol. I (Barcelona, 2000),  361--392, Progr. Math., 201, Birkhäuser, Basel, 2001.

%%\bibitem{gangl:GonMPMTM} {\bf Goncharov, A.B.}  {\it Multiple polylogarithms and mixed Tate motives.}  arXiv: math.AG/0103059 

%%\bibitem{gangl:GonFund} {\bf Goncharov, A.B.}  {\it Galois symmetries of fundamental groupoids and noncommutative geometry.} arXiv: math.AG/0208144 

\end{thebibliography}
\end{document}